\documentclass[11pt]{article}
\usepackage{srcltx}
\usepackage{times}
\usepackage{amssymb,amsfonts,amsmath,amsthm}
\theoremstyle{plain}
\newtheorem{thm}{Theorem}
\begin{document}
\title{\sc stable configurations of repelling points on flat tori}
\author{Marina Nechayeva and Burton Randol}
\date{}
\maketitle
\let\thefootnote\relax
\let\thefootnote\relax
\abstract{Flat tori are analyzed in the context of an intrinsic Fourier-analytic approach to electrostatics on Riemannian manifolds, introduced by one of the authors in 1984 and previously developed for compact hyperbolic manifolds.}
\bibliographystyle{plain}

\section{Introduction}
This paper continues themes introduced in \cite{randolchapter} and \cite{randol2014a}, whose notation we employ. We begin by providing proofs as well as extending to the flat case results from these references, which were developed in a hyperbolic context. In particular, we provide a complete proof of a general equilibrium condition which was stated but not proved in \cite{randolchapter}, and of which only a very brief sketch for hyperbolic manifolds was provided in \cite{randol2014a}. We then investigate several applications of these ideas to the case of flat tori, beginning with examples,  possibilities, and equidistribution results for the 1-torus, and continuing with a  discussion of the general case of the n-torus, which introduces new features. We conclude with an examination of questions about equidistribution, motivated by physics, that arose in the course of discussions with Dennis Sullivan, involving possibilities for extending our results to singular force laws, for example the Coulomb law on $T^3$, for which we develop an approximative approach.

Recall that \cite{randolchapter} introduced an intrinsic interpretation of electrostatics for certain types of Riemannian manifolds, in which force between two interacting points propagates, in a vectorially additive way, along the countably numerous geodesic segments connecting the two points, with magnitudes determined by a force law that is a function of distance. This differs from classical electrostatic settings, e.g., the Thomson problem, in which the manifold is isometrically embedded, generally in $R^3$, and the force between two points propagates along the unique single geodesic segment in $R^3$ connecting the points. 

We also remark at the outset that it is helpful to clearly differentiate the electrostatic model in which we are working from other intrinsic approches to electrostatics in the literature, in particular, from the very different shortest geodesic energy functional method, since the two are occasionally confused. To clarify the substantial differences between the two approaches, we note that instead of beginning with an arbitrary given energy functional, our starting point assumes a potential and its associated force law, which then produces an energy functional \cite{randolchapter} as a derivative secondary effect of the force law, and most importantly, in determining interaction between points, our method counts the effects that arise from all geodesics connecting them, which are ordinarily infinite in number. It is generally adapted to homogeneous spaces with a well-developed harmonic analysis, for which any two points in the universal covering space have a unique connecting geodesic, and it accommodates a phenomenon of non-trivial self-action, in particular in the hyperbolic case (cf.\ \cite{randol2014a}), although this phenomenon does not occur in the flat case, nor for any cases in the shortest geodesic approach. The latter takes as its starting point a postulated energy functional, defined as the sum of a specified real-valued function $f$ evaluated over the finitely numerous globally shortest distances between distinct points of a configuration. This is analogous to the classical treatment in Euclidean space, where the shortest connecting geodesic in the containing space is the unique connecting geodesic, and thus the only one requiring consideration. Although this approach has many interesting and important applications, we note that in the particular context in which we are interested, i.e., intrinsically in a generally non-simply connected manifold, restricting the count to shortest geodesics can lead to problems if one wishes to define a naturally associated force field, since discontinuities and instabilities in the dynamics can occur when one reaches the cut locus around a point, even in the case of the 1-torus, and the approach in \cite{randol2014a} was, among other considerations, specifically designed to eliminate this issue for such cases. In our view, both approaches have appropriate venues and applications, and broadly speaking, the shortest distance approach employs convexity arguments, whereas the method that counts all geodesics is fundamentally Fourier-analytic in character.

As mentioned above, in some cases it is important in our intrinsic treatment to consider the way in which a point may act on itself, and in the usual embedded version of this type of question, in which the set of standardly considered interacting point-pairs consists of the Cartesian product, minus the diagonal, of the specified point ensemble with itself, the issue by definition does not arise. 

On the other hand, in some geometries, our consideration of the intrinsic problem naturally leads to the possibility of non-trivial self action (c.f. \cite{randol2014a}, p. 2772, remark 4). This notion requires a suitable definition and analysis, which is not present in the existing literature, apart from the brief mention of such a possibility in Remark 4 of \cite{randol2014a}. Because of this gap, and to give a more complete picture of the scope and applicability of our method, we include a short analysis of this phenomenon here, even though non-trivial self action does not occur in the flat case, since in those cases where it does occur, the locations of an energy minimizing self-acting point encode geometric information about the containing manifold and constitute an important aspect of the general theory. More generally,  the study of energy minimizing configurations containing an arbitrary number of points is of interest in all cases, including the flat case. It is important to note that such configurations depend on the particular underlying potential function  (cf.\ Theorems 1 and 6 of this paper).

\begin{sloppypar}In \cite{randolchapter} and \cite{randol2014a}, there is postulated an underlying function \mbox{$k(\rho) \; \rho \in (-\infty,\infty)$,} corresponding to a classical potential function, which in turn produces a force law $H(\rho)$, given by the relation $H(\rho) = -k'(\rho)$. The function $k(\rho)$ can be quite general, although it is subject to certain requirements which are intended to secure the absolute uniform convergence of series which arise in the theory, as well as the validity of several derived identities. The contextual framework of the discussion in \cite{randolchapter} was that of compact hyperbolic manifolds, and the underlying analytic tool was the Selberg pre-trace formula. In the flat torus case, which was not explicitly treated in \cite{randolchapter} and \cite{randol2014a}, but with which we will be predominantly concerned in this paper, the corresponding analytic tool is the Poisson summation formula.\end{sloppypar}

Our operative analytical techniques require smoothness assumptions on the function $k(\rho)$ that are fairly strict in the hyperbolic case and less so in the case of flat tori. Another feature of our treatment of the intrinsic case is that, as is usual and useful in many applications of the Selberg pre-trace and Poisson summation formulas, great variety is permitted in the selection of the function $k(\rho)$. We note that in some, but by no means all, studies of electrostatics, e.g., in the case of the Coulomb potential, the function $k(\rho)$ has a singularity at the origin, which is not initially covered by our treatment, and apart from the last section of this paper, we will restrict our discussion to the non-singular case.

Retaining for now an assumption in \cite{randolchapter} that $k(\rho)$  is the restriction to $[0,\infty)$ of an even $C^\infty$ function on $(-\infty,\infty)$, and hence that the force function $H(\rho) = -k'(\rho)$ vanishes at $\rho = 0$, 
we also make a temporary assumption that $x \neq y$. Then, as noted in \cite{randolchapter} for the hyperbolic case, under suitable additional assumptions on $k(\rho)$, the vectorial effect at $x$ of a unit charge located at $y$ is given by \begin{equation}\begin{split}\sum_\gamma H(d(x,\gamma y))\vec{V}_\gamma  & =\\ -\nabla_x \sum_\gamma k(x,\gamma y) & =    - \nabla_x \sum_{n=1}^\infty h(r_n) \varphi_n(x)\overline{\varphi}_n(y)\,, \label{eq1}\end{split}\end{equation}where in the hyperbolic case as treated in \cite{randolchapter} and \cite{randol2014a}, $\gamma$ ranges over the elements of the lattice whose action on the appropriate hyperbolic covering space $H^n$ gives $M$, and $\vec{V}_\gamma$ is the tangent vector at $x$ of the unique directed geodesic segment in $H^n$ from $y$ to $x$. In this context, $h$ is the Selberg transform of $k$, the $r_n$'s are the standard parametrization in the Selberg theory of the Laplace eigenvalues on $M$, the $\varphi_n$'s are the corresponding Laplace eigenfunctions on $M$, and in the last sum, half the $r_n$'s are counted and the term corresponding to the constant eigenfunction is omitted, since the gradient operation removes the necessity for its presence.

In the case in which $M$ is a flat torus of unit volume given by the action of a lattice $\Gamma$ on a Euclidean space, an argument effectively identical to that of \cite{randolchapter} using the Poisson summation formula shows that an analog of formula (\ref{eq1}) is valid, in which $h$ is the Fourier transform of the radial function $k$ and the sum is taken over the non-zero elements $\nu$ of the dual lattice of $\Gamma$, with $r_n$ replaced by $\nu$ and $\varphi_\nu$ by $e^{2\pi i(\nu ,x)}$.

Note that under the currently operative assumption $H(0)=0$, the formula in (\ref{eq1}) has an obvious natural meaning if $x=y$, since the term corresponding to $\gamma = \mathrm{identity}$, namely $H(0)\vec{V}_{\mathrm{identity}}$\,, vanishes, and is thus not present. Moreover, this condition is continuous, in the sense that as $y \rightarrow x$, the corresponding term tends to zero. Thus, the effect of a charge at $x$ on itself has a natural definition as 
\[\sideset{}{'}\sum_{\gamma} H(d(x,\gamma x))\vec{V}_\gamma\,,\] where the prime indicates that $\gamma = \mathrm{identity}$ is omitted from the sum, and by either the pre-trace or the Poisson formula, the correct definition on the eigenfunction side is \[\nabla_x \sum_{n=1}^\infty h(r_n)|\varphi_n(x)|^2\,.\](Self action does not occur in the case of a flat torus, since the contribution from a geodesic segment is canceled by that from its negative counterpart. This also follows from the last formula, since the eigenfunctions in the flat case have constant absolute value $1$. Non-trivial self action can, however, occur in the hyperbolic case.)

In view of the above, the effect at $x_j$ of points at $x_1,\ldots x_N$ is given by \[-\nabla_{x_j} \sum_{n=1}^\infty h(r_n)\varphi_n(x_j)(\overline{\varphi}_n(x_1)+\cdots +\overline{\varphi}_n(x_N))\,,\]so the effect is null at $x_j$ if and only if \[-\nabla_{x_j} \sum_{n=1}^\infty h(r_n)\varphi_n(x_j)(\overline{\varphi}_n(x_1)+\cdots +\overline{\varphi}_n(x_N)) =0 \,.\]

If the configuration is in equilibrium at each of the $x_j$'s, then by adding and collecting terms, we find that \[-\nabla_{x} \sum_{n=1}^\infty h(r_n)(\varphi_n(x_1)+ \cdots \varphi_n(x_N))(\overline{\varphi}_n(x_1)+\cdots +\overline{\varphi}_n(x_N))=0\,,\]or  \[-\nabla_{x} \sum_{n=1}^\infty h(r_n)|\varphi_n(x_1)+ \cdots \varphi_n(x_N)|^2 =0\,,\]where the gradient can be regarded as taken over the Cartesian product $M^N$, since the gradient of the metric product is the orthogonal direct sum of the gradients in the factors.

Thus, a necessary and sufficient condition for the equilibrium of a configuration $\{x_1\ldots x_N\}$ is that it be a critical point on $M^N$ of 
\begin{equation}\sum_{n=1}^\infty h(r_n)|\varphi_n(x_1)+ \cdots \varphi_n(x_N)|^2 \label{energy1}\,,\end{equation} and as a corollary, a single point $x$ is in equilibrium under self action if and only if it is a critical point of\begin{equation}\sum_{n=1}^\infty h(r_n)|\varphi_n(x)|^2\label{energy_singlepoint}\,.\end{equation}

In particular, the above discussion strongly suggests that the quantities given by (\ref{energy1}) and its specialization (\ref{energy_singlepoint}) provide appropriate definitions for an energy functional in our setting, and that if $h$ is non-negative (repelling configurations), those configurations for which it is globally minimized will be of exceptional interest, while if $h$ is non-positive (attracting configurations), configurations for which it is globally maximized will play this role. Note in particular that in the hyperbolic case, the determination of points satisfying (\ref{energy_singlepoint}) is non-obvious, while in the flat torus case, (\ref{energy_singlepoint}) has the same value for all points, since $|\varphi_n(x)|^2 \equiv 1$.

We conclude this section by noting that we can associate to a potential function what may be termed its dual, which is simply the negative of the original function. As motivation, suppose for simplicity, although this is not required, that the force function $H$ is negative in $(0,\infty)$, i.e., not of mixed sign. Then the points in the configuration are mutually repelling, and the situation is electrostatic in character. If we now consider the environment produced by the dual potential, the points attract, and the situation becomes gravitational in character. The force vectors at each point have the same magnitudes in both cases, but are reversed in direction, from which it is clear that an electrostatically stable configuration, i.e., one in which the force resolves to zero at each point, is also gravitationally stable. This can be expressed as a general law: a configuration that is stable for a potential is also stable for its dual, or less formally, an electrostatically stable configuration is also gravitationally stable for the dual law.

\bigskip

\section{The 1-torus}In this section, we will look into the application of these ideas to the case in which $M$ is the unit volume 1-torus $T^1$, i.e., the quotient of $R^1$ by the integer lattice, in which case $h = \hat{k}$. For the 1-torus, the eigenfunctions are $\{e^{2\pi in}\}$, $n=0,\pm 1,\pm 2,\dots$, so in this case (\ref{energy1}) becomes\begin{equation}2\sum_{n=1}^\infty \hat{k}(n)|e^{2\pi inx_1}+ \cdots + e^{2\pi inx_N}|^2 \label{energy2}\,,\end{equation} which we will call, as suggested above, the energy of the configuration, and when $\hat{k} >0$, points on $M^N$ at which (\ref{energy2}) attains its global minimum $\mu_N$ are of exceptional interest.

 Note that the factor of $2$ in the the above sum is present because the absolute values in (\ref{energy2}) corresponding to $\pm{n}$ are identical, and the Fourier transform is defined with the scaling of $2\pi$ in the character, which in one dimension gives \[\hat{f}(\rho) = \int_{-\infty}^{\infty} f(t)e^{-2\pi i\rho t}\,dt\,.\]

\begin{thm}Suppose $\hat{k} >0$, and $\{S_N\}$ is a sequence of point configurations on $T^1$ which, for each $N$,
 globally minimize energy. Then $\{S_N\}$ is equidistributed as $N\rightarrow \infty$, at a rate that can be estimated above as a function of $N$.

\end{thm}   

This follows immediately from the theorem in pages 2770-2771 of \cite{randol2014a}.                   

\bigskip

We can easily derive a universal upper bound for $\mu_N$.  Namely, if we integrate the right side of (\ref{energy2}) over $T^N$, the $N$-fold Cartesian product of the $1$-torus $T^1$, then by orthonormality, bearing in mind that $|z|^2=z\overline{z}$, we obtain $2N\sum_{n=1}^\infty\hat{k}(n)$. It thus follows from the mean value theorem for integrals that there is a point $(p_1,\ldots,p_N) \in T^N $ such that 
\begin{equation}2\sum_{n=1}^\infty \hat{k}(n)|e^{2\pi inp_1}+ \cdots + e^{2\pi inp_N}|^2 = 2N\sum_{n=1}^\infty\hat{k}(n)\label{minimum}\,,\end{equation}so the right side of (\ref{minimum}) is an upper bound for $\mu_N$.  It is also obviously a lower bound for the global maximum of the energy achievable by any $N$-point configuration on $T^1$, although this is of limited interest, since it is obvious from (\ref{energy2}) that if, for example, $\hat{k}(n) \geq 0$,  the maximum energy is achieved when the points all coincide. 

Let us next examine how all of the above considerations apply to the case of $N$ equispaced points on $T^1$. For this case, of course, equidistribution in the limit is trivially obvious. Moreover, since all applicable results are invariant under isometry, the points can without loss of generality be taken to be the $N$th roots of unity. We also note that an equispaced configuration is in electrostatic equilibrium, or stable for short, since, if we regard the circle as standardly embedded and bounding a round disk in $R^2$, an equispaced configuration is invariant under reflection through a diameter of the disk passing through a point of the configuration, so all forces impinging on a point from one side are canceled by identical forces impinging from the other side. This last observation does not however, address the question of whether equispaced configurations are globally energy minimizing, although going in the other direction, it is a consequence of our previous discussion that global minimization, or for that matter, local minimization, implies stability, but it is on global minimization in the context of our model, specifically in the equispaced case, that we will next focus.

Accordingly, suppose that in (\ref{energy2}), $x_j = j/N\; (j=1,\ldots,N)$. Then, since a power of an $N$th root of unity is again an $N$th root of unity, and since when $\omega^N = 1$ and $\omega \neq 1$, $\sum_{j=1}^N \omega^j =0$, we conclude that all terms of (\ref{energy2}) vanish, except those occurring when $n= N,2N,3N\ldots$, in which case the corresponding term is $2N^2\hat{k}(n)$. In particular, the energy in the equispaced case is given by \begin{equation}2N^2\sum_{j=1}^\infty \hat{k}(jN)\label{equispaced energy}\,.\end{equation}

There are several interesting observations that we can now draw from this, depending on the behavior of the function $k$. In discussing these, we will work interchangeably with either $k$ or $\hat{k}$, depending on convenience. 

Setting $\nu_N$ to be the value of the energy (\ref{equispaced energy}) in the equispaced case, we note that $\mu_N \leq \nu_N$, and also that under mild conditions on the rate of decrease of $\hat{k}$, and the assumption that $\hat{k} \geq 0$, $\nu_N$ is asymptotically zero, in the sense that $\nu_N \rightarrow 0$ as $N \rightarrow \infty$. For example, if $\hat{k}(n) \ll \epsilon(N)/n^2$, with $\epsilon(N) \rightarrow 0$ as $N \rightarrow \infty$, then \[2N^2\sum_{j=1}^\infty \hat{k}(jN) \ll \epsilon(N)\,.\]On the other hand, if $1/n^{\alpha} \leq \hat{k}(n)$ for $\alpha < 2$, then $2N^2\sum_{j=1}^\infty \hat{k}(jN)$ tends to infinity as $N \rightarrow \infty$. In the transitional case for which $\hat{k}(n) $ is bounded between positive multiples of $1/n^2$, $\nu_N$  fluctuates between positive values or tends to a positive limit.

We now address the question of whether or not the equispaced configuration is globally energy minimizing in our model, and begin by showing that we can, by adding appropriate functions concentrated around a large $N=N_0$ to a suitable choice of $\hat{k}$, easily produce a situation in which $\hat{k}>0$, but $\nu_{N_0}$ is larger than the upper bound for $\mu_{N_0}$ given by the analysis leading to (\ref{minimum}). In particular, even if $\hat{k}>0$, the equispaced configuration, although stable, does not always globally minimize energy, although the corresponding potential function $k$ may not satisfy the classical positivity and monotonicity conditions. In order to exhibit this phenomenon, suppose $\psi(t)\geq 0$ is a standard even $C^\infty$ bump function supported in $[-\frac12,\frac12]$, with integral $1$ and $\psi(0)=1$, and $\hat{k}(t)$ is defined to be \mbox{$e^{-\pi t^2} + \frac12(\psi(t-N_0) + \psi(t+N_0))$.}Then $\hat{k}(N_0) \approx 1$, and $k(\rho) = e^{-\pi \rho^2}+\hat{\psi}(\rho)\cos 2\pi N_0\rho$. Now the contribution to the energy expression (\ref{equispaced energy}) from $j=1$ is $2N_0^2\hat{k}(N_0)\approx 2N_0^2$, and it is easy to see, even bearing in mind that the above choice of $\hat{k}$ depends on the choice of $N_0$, that the bound for $\mu_{N_0}$ given by (\ref{minimum}) is $\ll N_0$, so for large $N_0$ the equispaced configuration is not globally energy minimizing, which proves:

\begin{thm}The Fourier transform $\hat{k}$ can be chosen so that $\hat{k} > 0$ and the equispaced configuration is not globally energy minimizing.\label{notminimizing}\end{thm}

By contrast, we have:

\begin{thm} $\hat{k}$\ can be chosen to be of compact support, non-negative and non-increasing. with $k$ positive and strictly decreasing. These conditions imply, among other things, that for $N$ large, the equispaced configuration has zero energy, and is therefore globally energy minimizing.\label{thm3}\end{thm}

To show this, consider the Paley-Wiener type functions \[ \Phi_m(y) = -\int_{-\infty}^y \frac{sin^m t}{t^{m-1}}\,dt\,,\]where $m$ is an even integer greater than 2, to guarantee absolute convergence. This type of function, which for the case $m=4$ was introduced in \cite{khurgin-yakovlev} and is mentioned on page 82 of Higgins's survey article on cardinal series \cite{higgins}, provides a family of even, positive, band-limited functions which are strictly decreasing on $[0,\infty)$, and have interesting properties when chosen for the potential function $k$ in the equispaced case.

\begin{sloppypar}We quickly verify the above assertions about these functions. The positivity, evenness, and strict monotone decrease on $[0,\infty)$ are obvious, and it remains to verify that the functions are of Paley-Wiener type, i.e. band limited, meaning that their Fourier transforms have compact support. This follows from the facts that $\sin^m t/t^{m-1}$ is an entire function of $t$, and that the integral defining the function $\Phi_m(y)$ can be replaced, up to an additive constant, by an integral from $0$ to $y$, which implies that $\Phi_m(y)$ extends to an entire function of $y$, which can be defined up to an additive constant by an integral over the segment in $C^1$ connecting the origin to $y$. By an obvious estimate, $\Phi_m(y)$, regarded as a function in the complex plane, is of exponential type and is clearly $\ll y^{-(m-2)}$ and therefore $L^2$ on the real line, which by a Paley-Wiener theorem shows that $\Phi_m(y)$ is band limited. More elementarily, we could argue by noting that $\sin t/t$ is up to a positive multiple the Fourier transform of the indicator function of a compact interval, so $(\sin t/t)^m$ is in an obvious way the Fourier transform of a specific compactly supported $m$-fold convolution $\varphi(\rho)$. Moreover, integrating by parts,\[\int_{-\infty}^{\infty} \Phi_m(t)\,e^{-2\pi i \rho t}\, dt = \frac{1}{2\pi i\rho}\int_{-\infty}^{\infty} \Phi_m'(t)\,e^{-2\pi i \rho t}\,dt\]\\\[=-\frac{1}{2\pi i\rho}\int_{-\infty}^{\infty} \frac{sin^m t}{t^{m-1}}\,e^{-2\pi i \rho t}\, dt\]\\ \[= -\frac{1}{4\pi^2 \rho}\frac{d}{d\rho}\int_{-\infty}^{\infty}
\left(\frac{\sin t}{t}\right)^m\,e^{-2\pi i \rho t}\, dt\]\[= \frac{1}{4\pi^2 \rho}\,\varphi'(-\rho)\,,\]\\ which, since $\varphi'$ is a derivative of a compactly supported function and thus compactly supported, again shows that $\Phi_m(y)$ is band limited. Also, by examining the graphical interpretation of $\varphi(\rho)$, which is an $\frac{m}{2}$-fold convolution of a triangle function, it is easy to verify from the above that the Fourier transform of $\Phi_m$ is non-negative and non-increasing on $[0,\infty)$.\end{sloppypar}

Now take the potential function $k$ to be one of the $\Phi_m$'s. Then as noted, $k$ is positive and decreasing, and the transform $\hat{k}$, which figures in the energy expression (\ref{minimum}), is non-negative, non-increasing, and of compact support. This implies, since for roots of unity, the sum defining the energy reduces to $2N^2\sum_{j=1}^\infty \hat{k}(jN)$, that for sufficiently large $N$, depending on $k$, the energy in the equispaced case is zero, and is of course globally minimizing, which establishes Theorem \ref{thm3}.

\bigskip

We conclude with an additional criterion for when the equispaced configuration is globally energy minimizing. In connection with this, we will establish an interesting Fourier-analytic connection with a convexity result of G\"{o}tz \cite{gotz}, which he employed in the analysis of a very different electrostatic model (cf.\ our earlier remarks). In the case of our model, our examples are linked to the question of the convexity of the function on $[0,\frac12]$ whose Fourier coefficients coincide with $\{\hat{k}(n)\}_{n=1}^\infty$, and we will exploit the fact that in the flat torus case, the eigenfunctions are characters of the covering space,

We begin by noting that the energy expression (\ref{energy1}) can be rewritten, for $(x_1,\ldots ,x_N)$ in general position, as \[\sum_{n=1}^\infty \hat{k}(n)(e^{2\pi ix_1}+ \cdots + e^{2\pi inx_N})((e^{-2\pi inx_1}+ \cdots + e^{-2\pi inx_N}) \,,\]\[= N\sum_{n=1}^\infty \hat{k}(n) + \sum_{\stackrel{1\leq j,m \leq N}{j \neq  m}}\sum_{n=1}^\infty \hat{k}(n)\cos 2\pi n(x_m  -x_j)\]\[= \sum_{\stackrel{1\leq j,m \leq N}{j \neq  m}}\sum_{n=1}^\infty \hat{k}(n)\left(\frac{1}{N-1} +\cos 2\pi n(x_m -x_j)\right),\] where the last expression does not depend on the choice of representatives of $x_j$ and $x_m$ in the covering space.

\begin{sloppypar}In particular, for fixed $N$ the above expression can be written in the form \[\sum_{\stackrel{1\leq j,m \leq N}{j \neq  m}} f(d(x_j,x_m))\,,\]where $d$ is the length of a shortest geodesic, or arc, between $x_j$ and $x_m$, and the function $f$ is given by \[f(t)=  \sum_{n=1}^\infty \hat{k}(n)\left(\frac{1}{N-1} +\cos 2\pi nt\right) = c_N +\sum_{n=1}^\infty \hat{k}(n)\cos 2\pi nt \,.\]\end{sloppypar}

At this point we recall Proposition 9 from \cite{gotz} (cf.\ also \cite{saff2}), which is applicable to this situation. Specifically, \cite{gotz} is concerned with the study of configurations on the circle which absolutely minimize a functional that is a finite sum of the form
\[\sum_{j\neq m} f(d(x_j,x_m))\,,\] where $d(x_j,x_m)$ is the length of a shortest arc connecting $x_j$ and $x_m$, and $f$ is a real valued function on $[0,\frac12]$. 

In greater detail, Proposition 9 from \cite{gotz} states that if $f$ is convex and non-increasing, then the above functional attains a global minimum on the equispaced configuration, and if moreover $f$ is strictly convex, then the equispaced configuration is the unique such configuration.

We can now summarize the preceding discussion in the following theorem:

\begin{thm} If \[c_N +\sum_{n=1}^\infty \hat{k}(n)\cos 2\pi nt\] is the Fourier series of a convex non-increasing function on $[0,\frac12]$, then the equispaced configuration gives a global energy minimum in our sense. Furthermore, if the convex function is strictly convex, then the equispaced configuration is the unique global minimum for the potential $k$. (The additive constant is immaterial to this result, which only depends on the Fourier coefficients for $n \geq 1$)\,. \label{mincrit}\end{thm}

\noindent{\sc Remark.} Theorem \ref{mincrit}, which involves a hypothesis on the Fourier transform of the potential function $k$, shows that the intrinsic problem can, in some geometries, be somewhat ``finitized,'' but it must be borne in mind that the function $f$ is not finitely defined, and that the above connection is not applicable to, e.g., the hyperbolic case, because in that case the eigenfunctions do not combine in the requisite way. Also, importantly, the use of the above-referenced Proposition 9 from \cite{gotz} is not applicable to higher dimensional flat tori, because of the presence of geometric information in addition to distances. Finally, as Theorem \ref{notminimizing} makes clear, the condition expressed in Theorem \ref{mincrit} is not necessary and sufficient, and is one of several  that produce minimality. Theorem \ref{mincrit} is not vacuous, inasmuch as it is not difficult to give examples of potentials for which its hypotheses are satisfied.

\bigskip

\section{The higher dimensional case}We can define a  mapping to an energy functional in the multidimensional flat case as well, although in dimensions greater than 1, the result from \cite{gotz} is no longer applicable. Namely, suppose $T^r$ is the $r$-dimensional integral torus, and $(x_1\ldots,x_N)$ is an $N$-tuple of points on $T^r$. Then by imitating the derivation in the $1$-dimensional case for the energy functional $f$, we are led to an expression of the form\[N\sum_{n\neq 0}\hat{k}(n) + \sum_{\stackrel{1\leq j,m \leq N}{j \neq m}}\sum_{n\neq 0}\hat{k}(n) e^{2\pi i(n,x_m-x_j)}\,,\]where $n$ is summed over the non-zero integer lattice points in $R^r$. This suggests, for $j\neq m$, defining vectors $v_{j,m} = x_m-x_j$. The preceding expression can then be written \[N\sum_{n\neq 0}\hat{k}(n) + \sum_{\stackrel{1\leq j,m \leq N}{j \neq m}}\sum_{n\neq 0}\hat{k}(n) e^{2\pi i(n,v_{j,m})}\,,\]which can be combined into\[\sum_{\stackrel{1\leq j,m \leq N}{j \neq m}}\sum_{n\neq 0}\hat{k}(n) (c_r(N)+ \cos 2\pi (n,v_{j,m}))\,,\]where \[c_r(N)= 1/(N^{r-1} -1)\,.\]Thus if we define, for $x\in T^r$,\[f(x) = \sum_{\stackrel{1\leq j,m \leq N}{j \neq m}}\sum_{n\neq 0}\hat{k}(n) (c_r(N)+ \cos 2\pi (n,x))\,,\]we obtain, in several, dimensions, the following result:

\begin{thm} With notation as above, our energy expression is, up to an additive constant, equal to\[\sum_{\stackrel{1\leq j,m \leq N}{j \neq m}}f(v_{j,m})\,,\]which, as in the $1$-dimensional case, now gives a mapping from $k$ to an energy functional $f$ on $T^r$.\end{thm}

Note that in dimensions greater than $1$, $v_{j,m}$, and therefore $f(v_{j,m})$, depend on geometric information in addition to the distance between $x_j$ and $x_m$, so that, as mentioned, the result from \cite{gotz} is not applicable. Nevertheless. the following theorem is valid.

\begin{thm}Suppose $\hat{k} >0$, and $\{S_N\}$ is a sequence of point configurations on $T^r$ which, for each $N$,
 globally minimize energy. Then $\{S_N\}$ is equidistributed as $N\rightarrow \infty$, at a rate that can be estimated above as a function of $N$.\end{thm}

As in the case of $T^1$, this follows immediately from the theorem in pages 2770-2771 of \cite{randol2014a}.      

\bigskip

\section{Remarks on the Coulomb potential}

\bigskip

We end with some remarks that are pertinent to our discussion of the multidimensional case.They arose from correspondence with Dennis Sullivan, who was interested in possible extensions of the scope of our intrinsic electrostatic model \cite{randolchapter}. For example, is there some reasonable interpretation that can be associated with the Coulomb potential in this model. More specifically, is there such an interpretation in the case of a flat 3-torus, and if so, are the resulting globally minimizing configurations equidistributed in the limit?

We will briefly discuss some issues which make the provision of a suitable interpretative framework for this question difficult, and then sketch one of many possible approximations suggestive of information that can be developed in this direction.

\begin{sloppypar}We begin our remarks by noting an obvious difficulty posed by the Coulomb potential $1/r$ and its associated inverse square force law. Namely, $1/r$ does not fall within the scope of the admissible functions to which the theory has been developed, and an attempt to nevertheless simply apply it to the case of a flat 3-torus immediately runs into convergence problems. The (distributional) Fourier transform of the Coulomb potential in three dimensions is $1/r^2$. This well-known result in physics is often established by first noting that the Fourier transform of a Yukawa potential $e^{-\eta r}/r$ ($\eta >0$) is $1/(r^2 + \eta^2)$, and then dissolving $\eta$ to $0$ (note that the $1/r$ singularity is integrable around the origin in $R^3$).\end{sloppypar}

These facts prevent the use in this case of the equilibrium criterion given by (\ref{energy1}), since the number of corresponding dual lattice points grows too rapidly to accommodate its applicability. Note also that an attempt to define the total force profile of a point configuration arising from the exact Coulomb law will also encounter convergence problems, because of the growth with increasing length of the number of connecting geodesics that must be counted.

In view of this, it is reasonable to ask  whether or not there are suitable admissible approximations to the Coulomb potential that might prove useful in certain circumstances. For example, by mollifying the Coulomb potential at the two locations at which it is problematic, namely, at the origin and at infinity, it may be possible to obtain insight into questions for which the detection of a particular tendency is significant, as such approximate potentials approach the Coulomb potential.

To examine this approach,  mollification at the origin and at infinity is required, and there are many ways of doing this.  For this purpose, we illustrate with a Yukawa potential mollified at the origin. Specifically, in the setting of $R^3$, suppose \mbox{$\varphi(\rho) \leq 1$} is a non-decreasing radial $C^{\infty}$  function, which is $\equiv 0$ for $0 \leq \rho \leq \epsilon/2$, and $\equiv 1$ for $\rho \geq \epsilon$ ($\epsilon > 0$). Now consider the radial function $\varphi(\rho)\rho^{-1}e^{-\eta \rho}$ ($\eta >0$), and its Fourier transform \begin{equation}\int_{R^3}\varphi(\rho)\rho^{-1} e^{-\eta \rho}e^{-i(x,y)} dx\,,\end{equation}where for notational simplicity,  we are defining the Fourier transform without the $2\pi i$ factor in the exponent.

This is a Fourier-Bessel transform, since it is the Fourier transform of a radial function and hence itself radial, and in $R^3$, the Fourier-Bessel transform can be elementarily expressed, in this case as\[\frac{4\pi}{r}\int_0^\infty \varphi(\rho)e^{-\eta \rho}\sin(r\rho)d\rho\,.\]

This is a sine transform of a function on $[0,\infty)$. It behaves nicely at the origin and at infinity, and is an admissible $h$  in the formula (\ref{energy1}), so the previous techniques are applicable.

We are interested in the behavior of point configurations under the influence of these mollified potentials for small values of the parameters $\epsilon$ and $\eta$, and will call such potentials mollified Coulomb potentials, or just mollified potentials for short. In particular, we would like, if possible, to say something about whether asymptotic equidistribution occurs with such mollified potentials, as the number $N$ of points becomes large, subject to an appropriate minimization condition, perhaps one suggested by our previous study of configurations that globally minimize (\ref{energy1}).

\begin{sloppypar}There is a previously ultilized criterion for equidistribution that can be also adapted to the study of this question, in the form of the theorem contained in pages 2770-2771 of \cite{randol2014a}. An examination of its proof makes it clear that positivity plays a crucial role, which in the present context translates into a requirement that the Fourier transform $h(r)$ of a potential being investigated for an implication of equidistribution be positive, or more precisely, positive at the norms of the non-zero lattice points of $Z^3$.\end{sloppypar}

Rather than searching for a family of mollifications of the Coulomb potential having positive Fourier transforms, we will, in this brief discussion, adopt an approximate approach. Namely, the Weyl criterion shows that a sequence $\{S_1,S_2,\ldots\}$ of finite point sets of $M$ with cardinalities satisfying $|S_n| \uparrow \infty$, is equidistributed if, for all $j=1,2,\ldots$, it asymptotically integrates, with measure identically $1/N$ at each point of $S_N$, the Laplace eigenfunctions of $M$, ordered as in (\ref{eq1}), up to level $j$. With this in mind, we will call $\{S_1,S_2,\ldots\}$ approximately equidistributed to level $j_0$, if it asymptotically integrates all eigenfunctions up to and including level $j_0$.

We now take the mollified potentials introduced above,  and use, as counterparts and suitable approximations to the full energy functional (\ref{energy1}), its finite segments \begin{equation}\sum_{n=1}^{\ell} h_{\epsilon,\eta}(r_n)|\varphi_n(x_1)+ \cdots \varphi_n(x_N)|^2 \label{energy3}\end{equation}($\ell = 1,2\ldots$). 

As already noted, the Yukawa potential is in $L^1(R^3)$, and its Fourier transform in $R^3$ is  $1/(r^2 + \eta^2)$, which is positive on $[0,\infty)$. From this it is easy to show, for fixed $\eta > 0$, that there is an $r(\epsilon)$, with $r(\epsilon) \rightarrow \infty$ as $\epsilon \rightarrow 0$, such that  $h_{\epsilon,\eta}(r)$ is positive on $[0, r(\epsilon)]$. Keeping $\eta$ fixed, the following result now follows from an almost word-for-word adaptation of the proof of the theorem from \cite{randol2014a}, applied to a given finite segment $\mathcal{H}$ of the form (\ref{energy3}). 

\begin{thm}for any $\ell$ there is an $\epsilon(\ell) >0$ such that for any mollified potential with $0< \epsilon\leq\epsilon(\ell)$, a sequence $\{S_N\}$ of globally $\mathcal{H}$-minimizing configurations is approximately equidistributed to level $\ell$. \end{thm}

For purposes of this brief illustrative sketch, we will not here go into the further possibility of more intricate manipulations of $\epsilon$ and $\eta$, or into the investigation of other approximation methods, but will conclude by noting that for small $\eta$, the forces impinging on points of $S_N$ that arise from connecting geodesics with lengths in $[\epsilon,r(\epsilon)]$ closely approximate those that would be produced by the exact Coulomb force law.

\newcommand{\noopsort}[1]{}

\begin{flushleft}

Marina Nechayeva\\
LaGuardia Community College, NYC

\medskip

Burton Randol\\
CUNY Graduate Center, NYC

\end{flushleft}

\end{document}